\documentclass[a4paper,12pt]{article}
\usepackage{geometry,latexsym,amssymb,amsmath,color,graphicx}
\usepackage[latin5]{inputenc}
\usepackage{enumerate}
\usepackage{enumitem}
\usepackage[T1]{fontenc}
\usepackage{authblk}
\usepackage[all]{xy}
\usepackage{palatino}
\usepackage{indentfirst}
\usepackage{setspace}
\newcommand{\C}{\mathsf{C}}
\newcommand{\St}{\mathsf{St}}
\newcommand{\Cost}{\mathsf{Cost}}
\newcommand{\Ker}{\mathsf{Ker}}
\renewcommand{\Im}{\mathsf{Im}}
\newcommand{\CMod}{\mathsf{XMod}}
\newcommand{\XMod}{\mathsf{XMod}}

\newcommand{\Gpd}{\mathsf{GrpGpd}}
\newcommand{\Cat}{\mathsf{Cat}}
\newcommand{\E}{\mathsf{E}}
\newcommand{\TC}{\mathsf{Top^{\C}}}
\newcommand{\Cov}{\mathsf{Cov}}
\newcommand{\q}{\quad}
\newcommand{\NSGd}{\operatorname{NSGd}}
\newcommand{\NSCM}{\operatorname{NSCM}}

\newcommand{\Top}{\mathsf{Top}}
\newcommand{\Grp}{\mathsf{Grp}}

\newtheorem{example}{Example}[section]
\newtheorem{Def}[example]{Definition}
\newtheorem{Exam}[example]{Example}
\newtheorem{Prop}[example]{Proposition}
\newtheorem{Theo}[example]{Theorem}
\newtheorem{Lem}[example]{Lemma}
\newtheorem{Rem}[example]{Remark}
\newtheorem{Cor}[example]{Corollary}
\newenvironment{Prf}{{\bf Proof:} } {\hfill $\blacksquare$
	\mbox{}}
\begin{document}

\title{Normality and quotient in crossed modules, cat$^1$-groups and  internal groupoids within groups with operations}

\author[a]{Tunçar Şahan\thanks{E-mail : tuncarsahan@aksaray.edu.tr}}
\author[b]{Osman Mucuk\thanks{E-mail : mucuk@erciyes.edu.tr}}
\affil[a]{\small{Department of Mathematics, Aksaray University, Aksaray, TURKEY}}
\affil[b]{\small{Department of Mathematics, Erciyes University, Kayseri, TURKEY}}

\date{}

\maketitle
\begin{abstract} 
	In this paper we define the notions of  normal subcrossed module and quotient crossed module within groups with operations; and using the equivalence of   crossed modules over groups with operations and  internal groupoids we prove how normality and quotient concepts are related in these two categories. Further we prove an equivalence  of crossed modules over groups with operations  and cat$^1$-groups with operations for a certain algebraic category; and then by this equivalence we determine normal and quotient objects in the category of cat$^{1}$-groups with operations. Finally we characterize the coverings of cat$^{1}$-groups with operations.
\end{abstract}

\\ {\bf Classification:} 20L05, 57M10, 22AXX, 22A22 

\section{Introduction}

Crossed modules as   defined by  Whitehead \cite{Wth1,Wth2}
have been  widely used  in homotopy theory \cite{Brown-grenoble},
the theory of group representation (see \cite{BrownHubesshman} for a
survey),  in algebraic K-theory \cite{Loday}, and homological
algebra \cite{Hubeshman,Lue}. Crossed modules can be viewed as
2-dimensional groups \cite{BrLowDim}.

 The notions of subcrossed
module and normal subcrossed module were defined in \cite{KNor}. In \cite{BS} Brown and Spencer proved that the category of
internal groupoids  within the groups, which are also called in \cite{BS} under the name of  {\em $\mathcal{G}$-groupoids} and alternative names, quite generally used  are {\em group-groupoid} \cite{Br-Mu1}  or  {\em $2$-group} (see for example \cite{baez-lauda-2-groups})  is equivalent to the
category of crossed modules of groups.    Using the equivalence in  \cite{BS}, recently in \cite{Mu-Sa-Al}  normal and quotient objects in the category of group-groupoids  have been obtained.

In   \cite{Loday82}  Loday defined an algebraic object called {\em cat$^1$-group} as a group $G$ with two endomorphisms $s,t$ of $G$ such that
 $st=t$, $ts=s$ and  $[\Ker s, \Ker t]=0$, where $[\Ker s, \Ker t]$ represents the commutator subgroup of $G$; and proved that the categories  of cat$^1$-groups and crossed modules are equivalent.

 In \cite{Por} Porter proved a similar result to one in  \cite{BS}  holds for  certain algebraic categories, introduced by Orzech \cite{Orz}, which definition was adapted by him and called category of groups with operations. Applying Porter's result, the study of internal category theory was continued in the works of Datuashvili \cite{Kanex} and \cite{Wh}. Moreover, she developed cohomology theory of internal categories, equivalently, crossed modules, in categories of groups with operations \cite{Dat} and \cite{Coh}. The  equivalences of the categories  in \cite{BS} and  \cite{Por} enable  us to generalize some results on group-groupoids which are internal categories within groups to the more general internal groupoids for a certain  algebraic category $\C$ (see for example \cite{Ak-Al-Mu-Sa}, \cite{Mu-Be-Tu-Na},  \cite{Mu-Tu}  and \cite{Mu-Ak}).

 In this paper for an algebraic category $\C$ we define normal subcrossed module
and quotient crossed module for groups with operations; and then obtain  normal subgroupoid and  quotient  groupoid in internal groupoids corresponding
respectively to a normal subcrossed module and a quotient crossed module of groups with operations. Further we prove an equivalence  of crossed modules over groups with operations  and cat$^1$-groups with operations for a certain algebraic category. This equivalence enables us to determine normal and quotient objects; and coverings  in the category of cat$^{1}$-groups with operations.

Main results of this  paper  constitute some parts of the PhD thesis of first author at Erciyes University in 2014.

\section{Preliminaries}
Let  $G$  be a groupoid.  We write  $G_0$  for the set of
objects  of    $G $ and write $G_1$ for the set of morphisms. We also identify  $G_0$  with the set of
identities of  $G $ and so an  element  of  $G_0$  may be written as
$x$  or  $1_x$  as convenient.  We  write   $d_0, d_1 \colon
G_1\rightarrow G_0$  for the source and target maps, and, as usual,
write $G(x,y)$ for $d_0^{-1}(x)\cap d_1 ^{-1}(y)$, for $x,y\in G_0$.
The composition  $h\circ g$  of two elements of  $G$  is defined if
and only if  $d_0(h) =d_1(g)$, and so the   map $(h,g)\mapsto h\circ g$
is defined on the pullback  $G_1 {_{d_0}\times_{d_1}} G_1$
of $d_0$  and $d_1 $.  The {\em inverse} of $g\in G(x,y)$ is denoted by
$g^{-1}\in G(y,x)$.

If  $x\in G_0 $, we write $\St_Gx$  for $d_0^{-1}(x) $  and  call
the {\em star} of $G$ at $x$. Similarly  we write $\Cost_Gx$ for
$d_1^{-1}(x)$ and call {\em costar} of $G$ at $x$. The set of all
morphisms from $x$ to $x$ is a group, called {\em object group}
at $x$, and denoted by $G(x)$.

A groupoid $G$ is {\em totally intransitive} if $G(x,y)=\emptyset$
for all $x,y\in G_0$ such that $x\neq y$. Such a groupoid is determined entirely by the family
$\{G(x)\mid  x\in G_0\}$ of  groups. This totally intransitive
groupoid is sometimes called {\em totally disconnected} or {\em
bundle of groups} (Brown \cite[pp.218]{Br1}).

Let $G$ be a groupoid. A {\em subgroupoid} $H$ of $G$ is a pair of
subsets $H_1\subseteq G_1$ and $H_0\subseteq G_0$   such that
$d_0(H_1)\subseteq H_0$, $d_1(H_1)\subseteq H_0$, $1_x\in H_1$ for each
$x\in H_0$ and $H_1$ is closed under the partial multiplication and
the inversion in $G$. A subgroupoid $H$ of $G$ is called {\em wide}
if $H_0=G_0$.

\begin{Def}\label{normalsubgrpd}{\em Let $G$ be a groupoid.  A subgroupoid $N$
of $G$ is called {\em normal} if it is  wide in $G$ and   $g\circ N(x)=
N(y)\circ g$ for objects $x,y\in G_0$  and $g\in G(x,y)$.  } 
\end{Def}

%The normality condition can   also be  stated  as follows.
%\begin{Prop}\label{Propnormality}   Let $G$ be a groupoid and
%$N$ a subgroupoid of  $G$. Then $N$ is normal
%if and only if  for $a,b \in G(x,y)$ such that $a \circ b^{-1}\in N(x)$ and $c\in N(y)$, we have $a\circ c\circ b^{-1}\in
%N(x)$.\end{Prop}
%\begin{Prf} If  $N$ is normal,
%$a,b \in G(x,y)$ such that $ a \circ b^{-1}\in N(x)$  and $c\in
%N(y)$, then $a\circ c \circ a^{-1}\in  N(x)$. Since $N(x)$ is a
%group $a\circ c \circ a^{-1}\circ a\circ b^{-1}=a\circ c\circ b^{-1}\in  N(x)$. \end{Prf}

Quotient groupoid is  formed as follows (Higgins \cite[pp.86]{Hi} and Brown \cite[pp.420]{Br1}). Let  $N$ be a normal subgroupoid
of the groupoid  $G$.  The components of $N$ define a partition on
$G_0$ and we write $[x]$ for the class containing $x$.  Then $N$ also defines an
equivalence relation on $G_1$ by $g\sim h$ for $a,b \in G_1$ if and only
if $g=n\circ h \circ m$ for some $m,n\in N_1$.  A partial composition $[h]\circ[g]$ on the morphisms   is defined if and only if
there exist $g_1\in [g]$, $h_1\in [h]$ such that $h_1\circ g_1$ is
defined in $G_1$ and then $[h]\circ[g]=[h_1\circ g_1]$. This partial composition defines a groupoid on   classes $[x]$'s as objects.
The groupoid defined in this manner is called {\em quotient groupoid} and denoted by $G/N$.

As it is stated in Brown \cite[pp.218]{Br1},  in the case where the
normal subgroupoid $N$ is totally intransitive,  we have that
$(G/N)_0=G_0$ and  $G/N(x,y)$ consists of all cosets $g\circ N(x)$ for
$x,y\in G_0$ and $g\in G(x,y)$. The groupoid
composition becomes \[(h\circ N(y))\circ(g\circ N(x))=(h\circ g)\circ N(x).\]
for $g\in G(x,y)$ and  $h\in G(y,z)$.

We recall  that a crossed module of groups originally defined by  Whitehead  \cite{Wth1,Wth2},
 consists of two groups $A$ and $B$,
an action of $B$ on $A$ denoted by $b\cdot a$ for $a\in A$ and $b\in B$;
and a morphism  $\alpha\colon A\rightarrow B$ of groups satisfying the
following conditions for all $a,a_1\in A$ and $b\in B$
\begin{enumerate}[label=(\roman{*}), leftmargin=1cm]
\item $\alpha(b\cdot a)=b+\alpha(a)-b$,
\item $\alpha(a)\cdot a_1=a+a_1-a$.
\end{enumerate}
We will denote such a crossed module by $(A,B,\alpha)$. Let $(A,B,\alpha)$
and $(A',B',\alpha')$ be two crossed modules. A morphism  $(f_1,f_2)$
from $(A,B,\alpha)$ to $(A',B',\alpha')$  is a pair of morphisms of groups
$f_1\colon A\rightarrow A'$ and $f_2\colon B\rightarrow B'$  such that
$f_2\alpha=\alpha'f_1$ and  $f_1(b\cdot a)=f_2(b)\cdot f_1(a)$ for  $a\in A$ and $b\in B$.

It was proved by Brown and Spencer in \cite[Theorem 1]{BS} that the category $\CMod(\Grp)$ of crossed modules
over  groups is equivalent to the category $\Gpd$ of group-groupoids.

%By Loday in  \cite{Loday82} a new algebraic object called cat$^1$-group (or 1-cat-group)
%is defined as follows:
%A group $G$ is called a $cat^1-group$ if there are endomorphisms $s,t$ of $G$ such that
%\begin{enumerate}[label=(\roman{*}), leftmargin=1cm]
%\item $st=t$, $ts=s$;
%\item $[\Ker s, \Ker t]=0$.
%\end{enumerate}
%where $[\Ker s, \Ker t]$ represents the commutator subgroup of $G$. A cat$^1$-group
%morphism $f\colon G\rightarrow G'$ is a morphism of groups  which is compatible with
%the endomorphisms $s$ and $t$. The category of cat$^1$-groups is denoted by $\Cat^1(\Gr)$.
%In the same paper cited it was proved that the category $\Cat^1(\Gr)$ of cat$^1$-groups is
%equivalent to the category $\CMod(\Grp)$ of crossed modules of  groups.

\section{Normal and  quotient crossed modules in groups with operations}

The idea of the definition of categories of groups with
operations comes from Higgins \cite{Hig} and Orzech \cite {Orz};
and the definition below is from Porter \cite{Por} and Datuashvili
\cite[pp.~21]{Tamar}, which is adapted from Orzech \cite {Orz}.

\begin{Def}\label{groupswithoperations} {\em From now on  $\C$  will be a category of groups with a set  of  operations $\Omega$
and with a set $\E$  of identities such that $\E$ includes the group laws, and the
following conditions hold: If $\Omega_i$ is the set of $i$-ary operations in $\Omega$, then

(a) $\Omega=\Omega_0\cup\Omega_1\cup\Omega_2$;

(b) The group operations written additively $0,-$ and $+$ are  respectively
the  elements of $\Omega_0$, $\Omega_1$ and
$\Omega_2$. Let $\Omega_2'=\Omega_2\backslash \{+\}$,
$\Omega_1'=\Omega_1\backslash \{-\}$ and assume that if $\star\in
\Omega_2'$, then $\star^{\circ}$ defined by
$a\star^{\circ}b=b\star a$ is also in $\Omega_2'$. Also assume
that $\Omega_0=\{0\}$;

(c) For each   $\star \in \Omega_2'$, $\E$ includes the identity
$a\star (b+c)=a\star b+a\star c$;

(d) For each  $\omega\in \Omega_1'$ and $\star\in \Omega_2' $, $\E$
includes the identities  $\omega(a+b)=\omega(a)+\omega(b)$ and
$\omega(a)\star b=\omega(a\star b)$.

A category satisfying the conditions (a)-(d) is called a {\em category of groups with operations }.}\end{Def}

\begin{Rem}{\em
The set $\Omega_0$ contains exactly one element, the group identity;
hence for instance the category of associative rings with unit is not
a category of  groups with operations.}\end{Rem}

\begin{Exam}{\em
The categories of  groups, rings generally  without identity, $R$-modules,
associative, associative commutative, Lie, Leibniz, alternative algebras
are examples of categories of   groups with operations.}
\end{Exam}

A {\em  morphism} between any two objects of $\C$ is a group homomorphism, which preserves the operations in $\Omega_1'$ and $\Omega_2'$.

The topological version of this definition can be stated as follows:

\begin{Def}{\em Let $X$ be an object in $\C$. If $X$ has a topology such that all  operations in
$\Omega$ are continuous, then $X$ is called a {\em topological group
with operations} in $\C$.  }\end{Def}

In particular if $\C$ is the category of groups, then a topological group
with operations just becomes  a topological group and if  $\C$ is the category of  R-modules, then it becomes a topological R-module.

We will denote the category of topological groups with operations by $\TC$.

For  the  objects $A$ and  $B$  of $\C$,  the direct product  $A\times B$  with the usual operations becomes a group with operations.   Hence the category $\C$ has finite products.

The subobject in the category $\C$ can be defined as follows.
\begin{Def}{\em   Let $A$ be an object in $\C$. A subset   $B\subseteq A$
is called a {\em subgroup with operations of} $A$ if the following conditions are satisfied:
\begin{enumerate}[label=(\roman{*}), leftmargin=1cm]
  \item $b\star b_1\in B$ for  $b,b_1\in B$ and  $\star\in\Omega_2$;
  \item $\omega(b)\in B$ for $b\in B$ and  $\omega\in\Omega_1$.
\end{enumerate}}\end{Def}

The normal subobject in the category $\C$ of groups with operations is defined as follows.
\begin{Def} \emph{ \cite[Definition 1.7]{Orz} Let $A$ be an object in $\C$
and  $N$ a subgroup with operations of $A$.   $N$ is called a {\em normal
subgroup with operations} or  an {\em ideal} of $A$ and written  $N\lhd A$ if the following conditions are satisfied:
\begin{enumerate}[label=(\roman{*}), leftmargin=1cm]
  \item $(N,+)$ is a normal subgroup of $(A,+)$;
  \item  $a\star n\in N$ for  $a\in A$, $n\in N$ and $\star\in\Omega_2'$.\end{enumerate}}
\end{Def}

For a morphism  $f\colon A\rightarrow B$ in $\C$, $\Ker f=\{a\in A~|~f(a)=0\}$
is an ideal of  $A$.

A quotient object in $\C$ is constructed as follows:  Let $A$ be an
object in $\C$ and $N$ an ideal of $A$. Then the  relation on $A$
defined by  \[a\sim a_1 \q\text{iff}\q a-a_1\in N\]  is an equivalence
relation.  Then the quotient set   $A/N$  along with the
operations defined by
\begin{eqnarray*}
  [a]{\star}[a_1] &=& [a\star a_1] \\
 {\omega}([a]) &=& [\omega(a)]
\end{eqnarray*}
for  $\star\in\Omega_2$,  $\omega\in\Omega_1$ becomes an object in
$\C$ and called {\em quotient group with operations of} $A$ {\em by} $N$.

In the following proposition we prove that the category $\C$ has kernels in categorical sense.
\begin{Prop}
Let $A$  be an object in $\C$. Then  $N$ is an  ideal of $A$ if and only
if it is a kernel of a morphism in $\C$.
\end{Prop}

\begin{Prf}
We have already seen that the kernel of a morphism in $\C$ is an ideal of the domain.

Conversely if  $N$ is an ideal of $A$, then the quotient $A/N$ becomes an object in $\C$ and quotient morphism $p\colon A\rightarrow A/N$
has $N$ as kernel.
\end{Prf}

Let $A$ and $B$ be two groups with operations in $\C$. An {\em extension} of $A$ by $B$ is
an exact sequence
\begin{align}
0\longrightarrow A\stackrel{\imath}\longrightarrow E\stackrel{p}\longrightarrow B\longrightarrow 0
\end{align}
in which $p$ is surjective and $\imath$ is the kernel of $p$.  It is {\em split } if there is a morphism $s\colon  B \to E$ such
that $p s = \imath d _B$.  A split extension of $B$ by $A$ is called a {\em  $B$-structure} on $A$.  Given  such a $B$-structure on $A$ we  get actions of $B$ on $A$ corresponding to the operations in $\C$. For any $b\in B$, $a\in A$ and $\star\in \Omega'_2$ we have the actions called {\em derived actions} by Orzech \cite[pp.293]{Orz}
\begin{equation} \label{eq12} \begin{array}{rcl}
b\cdot a & = &s(b)+a-s(b)\\
 b\star a  & = &s(b)\star a.
   \end{array}\end{equation}

Given a set of actions of $B$ on $A$ (one for each operation in $\Omega_2$), let
$A\rtimes B$ be a universal algebra whose underlying set is $A\times B$ and whose
operations are
\begin{eqnarray*}
(a', b') + (a, b) &=& (a' + b' \cdot a,  ~b' + b), \\
(a', b') \star (a, b) &=& (a' \star a + a' \star b + b' \star a, ~b' \star b).
\end{eqnarray*}

\begin{Theo} \label{Derivedactionsemidrect}{\em \cite[Theorem 2.4]{Orz}} \label{Theoderivedaction}
A set of actions of $B$ on $A$ is a set of derived
actions if and only if $A\rtimes B$ is an object of $\C$.
\end{Theo}

We recall that for groups with operations  $A$ and $B$, in \cite[Proposition 1.1]{Dat} all necessary and sufficient conditions for the actions of $B$ on $A$ to be derived actions are determined.

\begin{Lem}\label{idealsemidrect}
Let $A,B,S$ and $T$ be objects in $\C$. Suppose that we have a set of derived actions of $B$ on $A$ and
a set of derived actions of $T$ on $S$. If $S\rtimes T$ is an ideal of $A\rtimes B$ then the followings
are satisfied.

\begin{enumerate}[label=(\alph{*}), leftmargin=1cm]
  \item  $S$ and $T$ are ideals of $A$ and $B$, respectively,
  \item  $b\cdot s\in S$ for all $b\in B$, $s\in S$,
  \item  $(t\cdot a)-a\in S$ for all $t\in T$, $a\in A$,
  \item  $b\star s\in S$ for all $b\in B$, $s\in S$,
  \item  $t\star a\in S$ for all $t\in T$, $a\in A$.
\end{enumerate}
\end{Lem}

\begin{Prf} In the  following proofs we assume that $S\rtimes T$ is an ideal of $A\rtimes B$ and $\star\in\Omega_2'$.
\begin{enumerate}[label=(\alph{*}), leftmargin=1cm]
  \item For $(s,0)\in S\rtimes T$ and $(a,0)\in A\rtimes B$, we have   \[(a,0)+(s,0)-(a,0)=(a+s-a,0)\in S\rtimes T\]
   and \[(a,0)\star(s,0)=(a\star s,0)\in S\rtimes T.\]
   Hence $a+s-a\in S$ and  $S$ is an ideal of $A$.
        Similarly if $(0,t)\in S\rtimes T$ and $(0,b)\in A\rtimes B$, then
        \[(0,b)+(0,t)-(0,b)=(0,b+t-b)\in S\rtimes T\] and \[(0,b)\star(0,t)=(0,b\star t)\in S\rtimes T.\]
        Hence $T$ is an ideal of $B$.
  \item If  $(s,0)\in S\rtimes T$ and $(0,b)\in A\rtimes B$, then
        \[(0,b)+(s,0)-(0,b)=(b\cdot s,0)\in S\rtimes T \]  and hence  $b\cdot s\in S$.
  \item  For  $(0,-t)\in S\rtimes T$ and $(t\cdot a,0)\in A\rtimes B$, it implies that
        \[(t\cdot a,0)+(0,-t)-(t\cdot a,0)=((t\cdot a)-a,-t)\in S\rtimes T\]  and so  $(t\cdot a)-a\in S$.
  \item For $(s,0)\in S\rtimes T$ and $(0,b)\in A\rtimes B$, we have
        \[(0,b)\star(s,0)=(b\star s,0)\in S\rtimes T \] and therefore $b\star s\in S$.
  \item If  $(0,t)\in S\rtimes T$ and $(a,0)\in A\rtimes B$, then         \[(0,t)\star(a,0)=(t\star a,0)\in S\rtimes T.\] and  $t\star a\in S$.
\end{enumerate}
\end{Prf}

Let $\mathcal{E}$ be a $B$-structure on $A$ and let $\mathcal{F}$ be a $T$-structure on $S$ as  below
\[\xymatrix{
\mathcal{E} : & 0 \ar[r] & A \ar@{->}[r]^{\iota} & E \ar@{->}[r]_{p} & B \ar@/_/[l]_{s} \ar[r] & 0 \\
\mathcal{F} : & 0 \ar[r] & S \ar@{->}[r]^{\iota'} & F \ar@{->}[r]_{p'} & T \ar@/_/[l]_{s'} \ar[r] & 0. }\]
If  $F$ is an ideal of $E$ then $\mathcal{F}$ is a normal substructure of $\mathcal{E}$. Then we can construct the quotient
split extension as follows.
\[\xymatrix{
\mathcal{E/F} : & 0 \ar[r] & A/S \ar@{->}[r]^{\iota^{\ast}} & E/F \ar@{->}[r]_{p^{\ast}} & B/T \ar@/_/[l]_{s^{\ast}} \ar[r] & 0 }\]
That is, $\mathcal{E/F}$ is a $B/T$-structure on $A/S$.

\begin{Def} {\em \cite{Por}} \label{Defcrosmod} {\em A crossed module in $\C$ is a triple $(A,B,\alpha)$, where $A$ and $B$ are the objects of $\C$, $B$
acts on $A$, i.e., we have a derived action in $\C$,  and $\alpha\colon A\rightarrow B$ is a
morphism in $\C$ with the conditions:

\begin{enumerate}[label={\roman{*})}, leftmargin=1.5cm]
  \item [CM1.] $\alpha(b \cdot a) = b + \alpha(a) - b$;
  \item [CM2.] $\alpha(a)\cdot a'=a+a'-a$;
  \item [CM3.] $\alpha(a)\star a'=a\star a'$;
  \item [CM4.]  $\alpha(b\star a)=b\star \alpha(a)$, $\alpha(a\star b)=\alpha(a)\star b$
\end{enumerate}
for any $b \in B$, $a, a'\in A$, and $\star\in\Omega_2'$.}\end{Def}

A morphism $(A,B,\alpha)\rightarrow(A',B',\alpha')$ between two crossed modules is a pair
$f\colon A\rightarrow A'$ and  $g\colon B\rightarrow B'$ of the morphisms in $\C$ such that
\begin{enumerate}[label={(\roman{*})}, leftmargin=1cm]
  \item $g \alpha(a)=\alpha' f(a)$,
  \item $f(b\cdot a)=g(b)\cdot f(a)$,
  \item $f(b\star a)=g(b)\star f(a)$
\end{enumerate}
for any $b\in B$, $a\in A$ and $\star\in\Omega_2'$.

\begin{Def}{\em \label{subcrossedmodule}
We call a crossed module  $(S,T,\sigma)$ in $\C$ as a {\em subcrossed module} of a  crossed
module $(A,B,\alpha)$ in $\C$ if
\begin{enumerate}[label={ \roman{*})}, leftmargin=1.7cm]
  \item [SCM1.]  $S$ is a subobject of $A$;  and  $T$ is a subobject of $B$;
  \item [SCM2.] $\sigma$ is the restriction of $\alpha$ to $S$;
  \item [SCM3.] the action of $T$ on $S$ is induced by the action of $B$
on $A$.
\end{enumerate}
}\end{Def}

\begin{Def}\label{nrmlsbcrsmd}{\em A subcrossed module $(S,T,\sigma)$ of $(A,B,\alpha)$ in the sense of Definition \ref{subcrossedmodule} is called {\em normal}
if
\begin{enumerate}[label=, leftmargin=1.5cm]
  \item [NCM1.]  $T$ is an ideal of $B$,
  \item [NCM2.]  $b\cdot s\in S$ for all $b\in B$, $s\in S$,
  \item [NCM3.]  $(t\cdot a)-a\in S$ for all $t\in T$, $a\in A$,
  \item [NCM4.]  $b\star s\in S$ for all $b\in B$, $s\in S$,
  \item [NCM5.]   $t\star a\in S$ for all $t\in T$, $a\in A$.
\end{enumerate}
}\end{Def}

\begin{Rem}{\em
Here we note that if $(S,T,\sigma)$ is a normal subcrossed module of $(A,B,\alpha)$
then by the conditions [CM2] of  Definition \ref{Defcrosmod};  and [NCM2] and [NCM4]
of  Definition \ref{nrmlsbcrsmd},  $S$ becomes an ideal of $A$.}
\end{Rem}

As an example if $(f,g)\colon(A,B,\alpha)\rightarrow (A',B',\alpha')$ is  a  morphism of crossed modules in $\C$, then $(\Ker f,\Ker g,\alpha_{|\Ker f})$, the kernel of $(f,g)$, is a normal subcrossed module of $(A,B,\alpha)$.

As a corollary of Lemma \ref{idealsemidrect} a  normal crossed module in $\C$ can be characterized as follow.

\begin{Cor}
A subcrossed module $(S,T,\sigma)$ of $(A,B,\alpha)$ is normal if and only if $S\rtimes T$ is an ideal of $A\rtimes B$.
\end{Cor}

We now obtain quotient crossed module in $\C$ as follows:
\begin{Theo} Let $(S,T,\sigma)$ be a normal subcrossed module of
$(A,B,\alpha)$ in $\C$. Then we have  a crossed module $(A/S,B/T,\alpha^{\ast})$ called
{\em quotient crossed module} where $A/S$ and $B/T$ are quotient groups with operations.
\end{Theo}
\begin{Prf}The actions of $B/T$ on $A/S$ are defined by
\begin{align*}
  [b]\cdot[a]& =  [b\cdot a] \\
  [b]\star[a]& =  [b\star a].
\end{align*}
These actions are well defined. If $[b]=[b_1]\in B/T$ and $[a]=[a_1]\in A/S$, then  $b-b_1\in T$ and $a-a_1\in S$.  Hence
\begin{align*}
  b\cdot a - b_1\cdot a_1 & = b\cdot(a-a_1+a_1) - b_1\cdot a_1 \\
                          & = b\cdot(a-a_1) + b\cdot a_1 - b_1\cdot a_1\\
                          & = b\cdot(a-a_1) + (b-b_1+b_1)\cdot a_1 - b_1\cdot a_1  \\
                          & = b\cdot(a-a_1) + (b-b_1)\cdot(b_1\cdot a_1) - b_1\cdot a_1
\end{align*}
and by the  substitutions $s=a-a_1\in S$, $t=b-b_1\in T$ and $a_2=b_1\cdot a_1\in A$,   we get that
\begin{align*}
  b\cdot a - b_1\cdot a_1 & = b\cdot(a-a_1) + (b-b_1)\cdot(b_1\cdot a_1) - b_1\cdot a_1 \\
                          & = b\cdot s + (t\cdot a_2) - a_2.
\end{align*}
Since $b\cdot s, (t\cdot a_2) - a_2 \in S$ by the conditions [NCM2] and [NCM3], we have that
$b\cdot s + (t\cdot a_2) - a_2=b\cdot a - b_1\cdot a_1 \in S$. This means $[b\cdot a]=[b_1\cdot a_1]$
in $B/T$. Moreover
\begin{align*}
  b\star a - b_1\star a_1 & = b\star a - b\star a_1 + b\star a_1 - b_1\star a_1 \\
                          & = b\star(a-a_1) + (b-b_1)\star a_1
\end{align*}
and by  the same substitutions above we get
\begin{align*}
  b\star a - b_1\star a_1 & = b\star s + t\star a_1.
\end{align*}
and so  by the conditions [NCM4] and [NCM5]  $b\star a - b_1\star a_1 \in S$.
Hence $[b\star a]=[b_1\star a_1]$ and therefore  these actions are well defined.

These  are derived actions  since the  conditions of   \cite[Proposition 1.1]{Dat} are satisfied.

On the other hand it is clear that the boundary map,
$\alpha^{\ast}\colon A/S\rightarrow B/T$ defined by $\alpha^{\ast}([a])=[\alpha(a)]$,
is well defined and the conditions [CM1]-[CM4] of  Definition \ref{Defcrosmod} are satisfied.
\end{Prf}

The  following result is useful for some proofs (see for example  the proofs of   Theorem \ref{corintgpd}  and  Theorem \ref{Thenormalsubcrosmod}).  The proof is clear and so it is omitted.
\begin{Prop} \label{exactseqofcrosmod} Let
\[
\xymatrix{
1 \ar[r] & S \ar[r] \ar[d]^{\sigma} & A \ar[r] \ar[d]^{\alpha} & P \ar[r]
\ar[d]^{\pi} & 1 \\
1 \ar[r] & T \ar[r] & B \ar[r] & Q \ar[r] & 1
}\]
be a short exact sequence of crossed modules in $\C$. Then $(S,T,\sigma)$ is a
normal subcrossed module of $(A,B,\alpha)$ and we have  the
short exact sequence of groups with operations
\[
\xymatrix{
1 \ar[r] & S\rtimes T \ar[r] & A\rtimes B \ar[r] & P\rtimes Q \ar[r] & 1
} \]
and so $S\rtimes T$ is an ideal of $A\rtimes B$.
\end{Prop}

In the following result we prove that a  normal subcrossed module in  $\C$ is categorically a  normal object in the category of crossed modules in $\C$.

\begin{Theo}
  $(S,T,\sigma)$  is a normal subcrossed module of $(A,B,\alpha)$  if and only if it is a kernel of a morphism $(f,g)\colon(A,B,\alpha)\rightarrow(C,D,\gamma)$
  of crossed modules in $\C$.
\end{Theo}

\begin{Prf}
We know that the kernel of a crossed module morphism $(f,g)\colon(A,B,\alpha)\rightarrow(C,D,\gamma)$ is a normal
subcrossed module of $(A,B,\alpha)$.

Conversely let $(S,T,\sigma)$ be a normal subcrossed module of $(A,B,\alpha)$.  Then  we have a morphism  $(p_A,p_B)\colon(A,B,\alpha)\rightarrow(A/S,B/T,\alpha^{\ast})$ of crossed modules whose kernel is the crossed module $(S,T,\sigma)$,  where  $p_A\colon A\rightarrow A/S$ and $p_B\colon B\rightarrow B/T$ are the natural projections.
\end{Prf}

\begin{Prop}\label{quotsemidirprod}
Let $(S,T,\sigma)$ be a normal subcrossed module of $(A,B,\alpha)$ in $\C$. Then the semi-direct groups with operations
$A/S\rtimes B/T$ and $ (A\rtimes B)/(S\rtimes T)$ are isomorph  in $\C$.
\end{Prop}

\begin{Prf}
It can be proved that  the function
\[\begin{array}{rcccc}
      \varphi & \colon & A/S\rtimes B/T & \longrightarrow & (A\rtimes B)/(S\rtimes T) \\
              &        & ([a],[b])      & \longmapsto     & [(a,b)]
    \end{array}\]
is an isomorphism in $\C$.
\end{Prf}

The isomorphism theorem for crossed modules in $\C$ can be given as follows.
\begin{Theo} \label{Isomthecrosmod}
Let $(f,g)\colon(A,B,\alpha)\rightarrow(A',B',\alpha')$ be a morphism of crossed modules; and let  $\Ker f=S$ and $\Ker g=T$. Then the image $(f(A),g(B),\alpha')$  is  a subcrossed module and isomorph to the
quotient crossed module $(A/S,B/T,\rho)$.
\end{Theo}

\begin{Prf}
It is easy to see that $(f(A),g(B),\alpha')$ is a subcrossed module and
\[(\widetilde{f},\widetilde{g})\colon (A/S,B/T,\rho) \rightarrow
(f(A),g(B),\alpha') \] is an isomorphism of crossed modules, where $\widetilde{f}(aS)=f(a)$ and
$\widetilde{g}(bT)=g(b)$ for $aS\in A/S$ and $bT\in B/T$.
\end{Prf}

\begin{Cor}
Let $(f,g)\colon(A,B,\alpha)\rightarrow(A',B',\alpha')$ be an epimorphism of crossed modules;
and let  $\Ker f=S$ and $\Ker g=T$. Then the image $(A',B',\alpha')$  is isomorph to the
quotient crossed module $(A/S,B/T,\rho)$.
\end{Cor}

\section{Fundamental crossed modules in  group with operations}

We recall a  major geometric example of a crossed module over groups  as follows: Let $(X,A,x)$  be a based pair of  spaces, where $X$ is a topological space and $x\in A\subseteq X$. As Whitehead proved the boundary map
\[\partial\colon \pi_2(X,A,x)\rightarrow \pi_1(A,x)\]
from the second relative  homotopy  group of $(X,A,x)$ to the fundamental group $\pi_1(A,x)$, together with the standard action of  $\pi_1(A,x)$ on $\pi_2(X,A,x)$ has the structure of crossed module. Here the elements of the second relative  homotopy  group $\pi_2(X,A,x)$ are the homotopy classes of the relative  paths as  pictured below and the compositions in the both directions  denoted induce the same group
\begin{figure}[ht] 
	\centering  \includegraphics{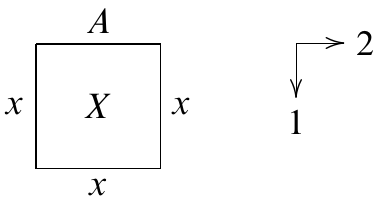}
\end{figure}
This is called {\em fundamental crossed module} of the based pair $(X,A,x)$ and denoted by $\Pi_2(X,A,x)$. In this manner we have a functor from  based pairs of topological  spaces to the crossed modules over groups
\[\Top_{\ast}^{2}\rightarrow \XMod (\Grp).\]
A 2-dimensional Seifert-van Kampen Theorem for fundamental crossed modules was proved in \cite[Theorem 2.3.1]{BHS}.

We know from \cite[Proposition 3.5]{Mu-Tu} that if $A$ is a topological group with operations in $\Top^{\C}$ , then  $\pi_1(A,0)$ is a group with operations in $\C$. We now prove that if $X$ is a topological group  with operations in $\Top^{\C}$  and $A$ is a subobject of $A$ in $\C$, then  $ \pi_2(X,A,0)$ is also a group with operations and the boundary morphism \[\partial\colon \pi_2(X,A,0)\rightarrow \pi_1(A,0)\] becomes a crossed module of groups with operations in $\C$.

As a preparation for the second relative paths $\alpha_1,\alpha_2$ and $\beta_1,\beta_2$ by the evaluation of the 2-dimensional paths
we obtain an interchange rule
\begin{align} \label{iterchange}
(\alpha_1 \circ \beta_1)+(\alpha_2\circ \beta_2)=(\alpha_1+\alpha_2)\circ (\beta_1+\beta_2)
\end{align}
whenever the compositions of the paths $\alpha_1 \circ \beta_1$ and $\alpha_2\circ \beta_2$ are defined.
Hence by the interchange rule (\ref{iterchange}) for second relative paths $\alpha, \beta$   we have that
\begin{align}
\alpha+ \beta\simeq (\alpha\circ \alpha_0)+(\alpha_0\circ \beta)=(\alpha+\alpha_0)\circ (\alpha_0+\beta)=\alpha\circ\beta
\end{align}
where $\alpha_0$ is the second zero path defined by $\alpha_0(s,t)=0$ for $0\leq s,t\leq 1 $. Hence two group operations on $\pi_2(X,A,0)$ are the same.
\begin{Theo} \label{Secondhomatpabel} If $X$ is a topological group with operations  in $\Top^{\C}$ and $A$ is  a   subobject   of $X$ in $\C$, then   $ \pi_2(X,A,0)$ is a group with operations in $\C$ and it is abelian with respect to ``$+$''. \end{Theo}
\begin{Prf} The binary operations on $\pi_2(X,A,0)$ are defined by $[\alpha]\star [\beta]=[\alpha \star \beta]$ for $\star\in \Omega_2$ and the unary operations are defined by $\omega [\alpha]=[\omega(\alpha)]$ for $\omega\in \Omega_1$. The other details are satisfied and hence $\pi_2(X,A,0)$ becomes a group with operations.

 Let $\alpha$ and $\beta$ be second relative paths. Then we define a homotopy
\begin{align}
F(r,s,t)=-\alpha(r,st)+\alpha(r,s)+\beta(r,s)+\alpha(r,st).
\end{align}
Here $F(r,s,0)= \alpha(r,s)+\beta(r,s)$ and  $F(r,s,1)=\beta(r,s)+ \alpha(r,s)$. Hence $\alpha+\beta$ and $\beta+\alpha$ are homotopic. Therefore  $\pi_2(X,A,0)$ is abelian with respect to ``$+$''.\end{Prf}

\begin{Theo} If $X$ is a topological group  with operations in $\Top^{\C}$  and $A$ is a subobject of $X$ in $\C$, then the boundary morphism \[\partial\colon \pi_2(X,A,0)\rightarrow \pi_1(A,0)\] becomes a crossed module of groups with operations in  $\C$.
\end{Theo}
\begin{Prf}  By Theorem  \ref{Secondhomatpabel},  $\pi_2(X,A,0)$ is abelian  with respect to ``$+$'' and for $[b]\in \pi_1(A,0)$,  $[\alpha]\in \pi_2(X,A,0)$  the actions  of $\pi_1(A,0)$ on $\pi_2(X,A,0)$  are defined by

\begin{align*}
  [b]\cdot[\alpha]& =  [\alpha] \\
  [b]\star[a]& = [\beta\star a]
\end{align*}
for $\star\in \Omega_2'$ where $\beta(s,t)=b(s)$ for $0\leq t \leq 1$. It is immediate to prove that these  are well defined derived actions and the conditions of Definition \ref{Defcrosmod} are satisfied. Therefore \[\partial\colon \pi_2(X,A,0)\rightarrow \pi_1(A,0)\]  becomes a crossed module of groups with operations.

\end{Prf}

Hence we have a functor from based pairs of topological groups with operations  to the crossed modules over groups with operations
\[(\Top^{\C})_{\ast}^{2}\rightarrow \XMod(\C)\]

\begin{Cor}
If $(X,A,0)$ is a based pair topological group with operations and  $(Y,B,0)$ is a normal subobject of $(X,A,0)$, then the crossed module $\Pi_2(Y,B,0)$ is a normal subcrossed module  of $\Pi_2(X,A,0)$ in $\C$.
\end{Cor}

\section{Normal  and quotient groupoids in the internal groupoids}
\begin{Def}\label{Internalgpd} {\em  An internal category $C$  in $\C$ is a category in which
the initial and final point maps $d_{0},d_{1}\colon C_1\rightrightarrows C_0$, the object inclusion map  $\epsilon\colon C_0\rightarrow C_1$
and the partial composition $\circ\colon C_1{_{d_0}\times_{d_1}}
C_1\longrightarrow C_1,(b,a)\mapsto b\circ a$ are the morphisms in the category
$\C$. }\end{Def} 

Note that since $\epsilon$ is a morphism in $\C$,
$\epsilon(0)=0$ and that the composition  $\circ$ being a morphism in $\C$,
implies that for all $a,b,c,d\in C$ and $\star\in \Omega_2$
\begin{align}\label{1}
(a\star b)\circ(c\star d)&=(a\circ c)\star(b\circ d)
\end{align}
whenever one side makes sense. This is called the {\em interchange law}.

As an application of the intercahange rule (\ref{1}) in an internal category $C$ for $a,b\in C_1$ such that $d_0(b)=d_1(a)=y$ we have the equality
\begin{equation}\label{comp}
\begin{array}{rl}
    b\circ a & = (b+1_0)\circ(1_y+(-1_y+a)) \\
             & = (b\circ 1_y) + (1_0\circ(-1_y+a)) \\
             & = b-1_y+a.
\end{array}
\end{equation}

As another easy application we note that any internal category $C$ in
$\C$ is an internal groupoid since given $a\in C$,
\begin{align}\label{2}
{a}^{-1}&=1_{d_1(a)}-a+1_{d_0(a)}\end{align}
satisfies $a^{-1}\circ a=1_{d_1(a)}$,
$a\circ a^{-1}=1_{ d_0(a)}$ and the map $C\rightarrow C,a\mapsto a^{-1}$ is also a morphism in $\C$.

In particular if $\C$ is the  category of groups, then an internal category $C$ in $\C$
becomes a group object in the category of groupoids, which is called as  {\em group-groupoid}, 2-{\em group} or $\mathcal{G}$-groupoids   \cite{BS}.

\begin{Exam} \label{fundgdgpgd}
{\em Let  $A$ be an object in $\C$. Then  the groupoid $G=A\times A$  is an
internal groupoid:  Here a pair $(a,b)$ is a morphism from $a$ to $b$  with
inverse morphism $(b,a)$. The  groupoid composition is defined by $(c,d)\circ (a,b)=(a,d)$
whenever $b=c$.  The induced  group operations are defined by
$(a,b)\star (c,d)=(a\star c,b\star d)$ for $\star\in \Omega_2$ and
$\omega(a,b)=(\omega(a),\omega(b))$ for $\omega\in \Omega_1$. }  
\end{Exam}

\begin{Exam} {\em  \cite[Example 3.8]{Ak-Al-Mu-Sa} If $X$ is  an object in $\TC$, then the fundamental groupoid $\pi X$ is an internal groupoid
in $\C$.}\end{Exam}

\begin{Rem} \label{Remcorfromdef}{\em We emphasize the  following points  from  Definition \ref{Internalgpd} \cite [Remark 3.7]{Ak-Al-Mu-Sa}:

\begin{enumerate}[label=(\roman{*}), leftmargin=1cm]
\item By Definition \ref{Internalgpd} we know  that in an internal groupoid $G$ in $\C$, the initial and final point maps $d_0$ and $d_1$, the object inclusion  map $\epsilon$ are the morphisms in $\C$ and the interchange law  (\ref{1}) is satisfied. Therefore in an internal groupoid $G$,  the unary operations are endomorphisms of the underlying groupoid of $G$ and the binary operations are  morphisms from the underlying groupoid of  $G\times G$ to the one of $G$.
\item Let $G$ be an internal groupoid in $\C$ and $0\in G_0$ the identity element. Then $\Ker d_0=St_G0$, called in \cite{Br1}  {\em transitivity component} or {\em connected component} of $0$, is also an internal groupoid which is also an ideal of $G$.
\end{enumerate}}
\end{Rem} 

\begin{Lem}\label{idealsemi}
Let $G$ be an internal groupoid in $\C$ and $N$ a wide subgroupoid of $G$. Then $N$ is a normal subgroupoid of $G$ in the sense of Definition \ref{normalsubgrpd} if and only if \[-1_x+N(x)=-1_y+N(y)\] for all
$x,y\in G_0$ with $G(x,y)\neq\emptyset$.
\end{Lem}

\begin{Prf}
Let $N$ be a normal subgroupoid of $G$. Then $g\circ N(x)=N(y)\circ g$ or by Equation  (\ref{comp}) equivalently
\[g-1_x+N(x)=N(y)-1_y+g\] for  $x,y\in G_0$ and $g\in G(x,y)$. Hence
\[ g+(-1_x+N(x))+(-g+1_y)  = N(y)\]
 and here all morphisms of $-1_x+N(x)$ are in $\Ker d_0$ and $-g+1_y \in\Ker d_1$. So  morphisms
of $-1_x+N(x)$ commute with $-g+1_y$.  Therefore we can write
    \begin{align*}
        g-1_x+N(x) & = N(y)-1_y+g \\
        g+(-g+1_y)+(-1_x+N(x)) & = N(y) \\
        1_y-1_x+N(x) & = N(y) \\
        -1_x+N(x) & = -1_y + N(y).
    \end{align*}

Conversely let \[-1_x+N(x)=-1_y+N(y)\] for all $x,y\in G_0$ with $G(x,y)\neq\emptyset$.
By reversing the steps above we get $g\circ N(x)=N(y)\circ g$ for objects $x,y\in G_0$ and $g\in G(x,y)$.
Hence $N$ becomes a normal subgroupoid of $G$.
\end{Prf}

The details of the following theorem is given in \cite{Akiz}.

\begin{Prop}{\em }
Let $G$  be an  internal groupoid in the groups with operations in $\C$  and $N$  a normal subgroupoid of $G$.  If  $N_1$  is a subobject of $G_1$ in $\C$,   then  the quotient groupoid $G/N$ becomes an internal groupoid in the groups with operations.
\end{Prop}

The following theorem was proved  in \cite[Theorem1]{Por}. Since some details of the proof will be used in later proofs we give a sketch proof.
\begin{Theo} \label{Theocatequivalencess}
The category $\CMod(\C)$ of crossed modules  and the category $\Cat(\C)$  of internal groupoids in $\C$ are equivalent.
\end{Theo}

\begin{Prf}
A functor $\delta\colon \Cat(\C)\rightarrow \CMod(\C)$ is defined as follows:
For an internal groupoid $G$,  let  $\delta(G)$ be the crossed module
$(A,B,d_1)$ in $\C$,  where $A=\Ker d_0$, $B=G_0$ and $d_1\colon A\rightarrow B$ is the
restriction of the target point map. Here $A$ and $B$ inherit the structures of group with
operations  from that of $G$, and the target point map
$d_1\colon A\rightarrow B$ is a morphism in $\C$. Further the  actions $B\times A\rightarrow A$ on
the group with operations $A$ are defined by
\begin{eqnarray*}
  b\cdot a &=& \epsilon(b)+a-\epsilon(b) \\
  b\star a &=& \epsilon(b)\star a
\end{eqnarray*}
for  $a\in A$, $b\in B$. The axioms of Definition \ref{Defcrosmod} are satisfied. Thus $(A,B,d_1)$ becomes a crossed module in $\C$.

Conversely define a functor $\eta\colon \CMod(\C)\rightarrow \Cat(\C)$ in the following way. For a crossed module  $(A,B,\alpha)$ in $\C$,  define an internal groupoid   $\eta(A,B,\alpha)$  whose  set of objects is the group with operations $B$ and  set of morphisms is the semi-direct product
$A\rtimes B$ which  is a group with operations by Theorem \ref{Derivedactionsemidrect}.
The  source and target  point maps are defined to be $d_0(a,b)= b$ and  $d_1(a,b)= \alpha(a)+b$  while the object inclusion map and groupoid
composition is given by $\epsilon(b)=(0,b)$ and \[(a_1,b_1)\circ (a,b)=(a_1+a,b)\] whenever
$b_1=\alpha(a)+b$. Hence $\eta(A,B,\alpha)$ is an internal groupoid. The other details of the proof  is obtained from that of \cite[Theorem 1]{Por}.
\end{Prf}

\begin{Lem} \label{lemsub}\label{lemnorm}  Let $G$ be an internal groupoid in $\C$  and  $H$ a subgroupoid of $G$.
\begin{enumerate}[label=(\alph{*}), leftmargin=1cm]
\item If $H_1$ is a subobject of $G_1$, then $H_0$ is also a subobject of $G_0$.
\item If $H_1$ is a an ideal of $G_1$ then $H_0$ is also an ideal of $G_0$.
\end{enumerate}
\end{Lem}

\begin{Prf}
\begin{enumerate}[label=(\alph{*}), leftmargin=1cm]
\item Let $x,y\in H_0$. Since $H$ is a subgroupoid of $G$   $1_x,1_y\in H_1$ and since $H_1$ is a subobject of $G_1$
   $1_x\star 1_y = 1_{x\star y}\in H_1$ and $\omega(1_x)=1_{\omega(x)}\in H_1$.  So
    $x\star y\in H_0$ and $\omega(x)\in H_0$. Hence  $H_0$ is a subobject of $G_0$.
\item Let $x\in G_0$ and $y\in H_0$. In this case $1_x\in G_1$ and $1_y\in H_1$.
    Since $H_1$ is an ideal of $G_1$ we have that $1_x+1_y-1_x=1_{x+y-x}\in H_1$ and $1_x\star 1_y = 1_{x\star y}\in H_1$
    and so $x+y-x\in H_0$, i.e. $(H_0,+)$ is a normal subgroup of $(G_0,+)$, and $x\star y\in H_0$.
    Hence  $H_0$ becomes an ideal of $G_0$.
\end{enumerate}
\end{Prf}

By Theorem \ref{Theocatequivalencess} and Lemma \ref{lemsub} we relate the subobjects in these categories as follows:

\begin{Theo}
Let $(S,T,\sigma)$ be a subcrossed module of a crossed module $(A,B,\alpha)$ in $\C$.  Suppose that  $H$ and $G$ are respectively
the internal groupoids in $\C$ corresponding to these crossed modules. Then $H_1$ is a subobject of $G_1$ and $H_0$ is a subobject of $G_0$.
\end{Theo}

\begin{Prf} By the detail of the proof of  Theorem \ref{Theocatequivalencess}, we know that $H_1=S\rtimes T$ and $G_1=A\rtimes B$.
It is clear that $S\rtimes T$ is a subobject of $A\rtimes B$ since  $(S,T,\sigma)$ is a subcrossed module of $(A,B,\alpha)$ . Hence  $H_1$ is a subobject of $G_1$ and by Lemma \ref{lemsub} $H_0$ is also a subobject of $G_0$.
\end{Prf}

Hence the notion of internal subgroupoid of an internal groupoid in $\C$ can be stated as follows:

\begin{Def}
\emph{ Let $G$ be an internal  groupoid in  groups with operations in $\C$ and  $H$ a subgroupoid of $G$ such that
$H_1$ is a subobject of $G_1$. Then   $H$ is called an {\em  internal subgroupoid} of $G$.}
\end{Def}

By Theorem \ref{Theocatequivalencess} and Lemma \ref{lemnorm} we relate the normal subobjects
in these categories as follows:

\begin{Theo} \label{corintgpd} Let $(S,T,\sigma)$ be a normal subcrossed module of a crossed module $(A,B,\alpha)$ in $\C$.
Suppose that $N$ and $G$ are respectively the internal groupoids in $\C$ corresponding to these crossed
modules. Then $H_1$ is an ideal of $G_1$ and $H_0$ is also an ideal of $G_0$.
\end{Theo}

\begin{Prf}
By the proof  of  Theorem \ref{Theocatequivalencess} $N_1=S\rtimes T$ and $G_1=A\rtimes B$.
By  Proposition \ref{exactseqofcrosmod}  $S\rtimes T$ is an ideal of $A\rtimes B$.
Therefore $H_1$ is an ideal of $G_1$ and by Lemma \ref{lemnorm} $H_0$ is an ideal of $G_0$.
\end{Prf}

Hence the notion of internal normal subgroupoid of an internal groupoid  in $\C$ can be stated as follows:

\begin{Def} \label{intnormsubgpd} {\em  Let $G$ be an internal  groupoid in groups with operations in $\C$ and  $N$ a subgroupoid of $G$ such that  $N_1$ is an ideal of $G_1$.
Then  $N$ is called  an {\em internal  normal subgroupoid} of $G$.}
\end{Def}

\begin{Exam} \label{Examsupgpgpd2}{\em These are some examples of internal normal subgroupoids:

\begin{enumerate}[label=(\roman{*}), leftmargin=1cm]
\item   Let $A$ be an object in $\C$ and $B$ an ideal of $A$. Then  $N=B\times B$ as defined  in Example \ref{fundgdgpgd} becomes
an internal normal subgroupoid of  $G=A\times A$.

\item  Let $G$ and $H$  be two internal groupoids in $\C$ and $f\colon G\rightarrow H$ a morphism of  internal groupoids.
Then the kernel \[\Ker f=\{a\in G_1 ~|~ f(a)=1_{0}\in H_1\}\] of $f$ is an internal normal subgroupoid of $G$.

\item If  $X$ is an object in $\TC$  and $Y$  an ideal of $X$ in $\C$, then  the fundamental
 groupoid $\pi Y$ is an internal  normal subgroupoid of $\pi X$.
 
\end{enumerate}}
\end{Exam}

\begin{Theo} \label{qugpdofgpgd}  Let  $G$ be an internal groupoid in $\C$ and  $N$
an internal  normal subgroupoid of $G$.  Then the crossed module
corresponding to $N$ is a normal subcrossed module of the one
corresponding to $G$.\end{Theo}

\begin{Prf}
By Definition \ref{intnormsubgpd} $N_1$ is an ideal of $G_1$. Let $(A,B,\alpha)$
and $(S,T,\sigma)$  be the corresponding crossed modules to $G$ and $N$ respectively.
So $A=\Ker d_0$, $B=G_0$, $S=\Ker d_{0}\cap N_1$ and $T=N_0$. To prove that $(S,T,\sigma)$ is a
subcrossed module of $(A,B,\alpha)$, we need to show that $(S,T,\sigma)$ satisfies
the conditions of Definition \ref{nrmlsbcrsmd}.
\begin{enumerate}[label={NCM\arabic{*}.}, leftmargin=2cm]
\item We know that by Lemma \ref{lemnorm} $T$ is an ideal of $B$.
\item Let $b\in B$ and $s\in S$. Then by the proof of Theorem \ref{Theocatequivalencess}
    $b \cdot s  = 1_b+s-1_b$ where $1_b,-1_b\in G_1$. Since $N_1$ is an ideal of $G_1$ we
    have that $1_b+s-1_b=b \cdot s\in N_1$ and
    \begin{align*}
        d_0(1_b+s-1_b) & = d_0(1_b)+d_0(s)-d_0(1_b) \\
                                  & = b+0-b\\
                                  & = 0.
    \end{align*}
   Hence $b \cdot s \in S$.
\item Let $t\in T$ and $a\in A$. Then
    \begin{align*}
        d_0((t\cdot a)-a) & = d_0(1_t+a-1_t-a) \\
                          & = d_0(1_t)+d_0(a)-d_0(1_t)-d_0(a)\\
                          & = t+0-t-0 \\
                          & = 0
    \end{align*}
    and so $(t\cdot a)-a\in A$. Moreover since $a\in G_1$, $1_t\in N_1$ and $N_1$ is an
    ideal of $G_1$ it follows that $a-1_t-a\in N_1$ and since $1_t\in N_1$ it implies that
    $1_t+a-1_t-a=(t\cdot a)-a\in N_1$. So $(t\cdot a)-a\in S$.
\item Let $b\in B$ and $s\in S$. Then $b\star s=1_b\star s\in N_1$ since $1_b\in G_1$, $s\in N_1$
    and $N_1$ is an ideal of $G_1$. Also
    \begin{align*}
        d_0(b\star s) & = d_0(1_b\star s) \\
                          & = d_0(1_b)\star d_0(s)\\
                          & = b\star 0 \\
                          & = 0
    \end{align*}
    and so $b\star s\in A=\Ker d_0$. Thus $b\star s\in S$.
\item Let $t\in T$ and $a\in A$. Then $t\star a=1_t\star a\in N_1$ since $1_t\in N_1$, $a\in G_1$ and
    $N_1$ is an ideal of $G_1$. Further
    \begin{align*}
        d_0(t\star a) & = d_0(1_t\star a) \\
                          & = d_0(1_t)\star d_0(a)\\
                          & = t\star 0 \\
                          & = 0
    \end{align*}
    and so $t\star a=1_t\star a\in A=\Ker d_0$. Thus $t\star a\in S$
\end{enumerate}
Therefore $(S,T,\sigma)$ is a normal subcrossed module of $(A,B,\alpha)$.
\end{Prf}

As a result of Theorem \ref{corintgpd} and Theorem \ref{qugpdofgpgd} we can state the following corollary.

\begin{Cor}
Let $G$ be an internal groupoid in $\C$ and $(A,B,\mu)$  the crossed module corresponding to $G$.
Then the category $\NSGd_{\Cat(\C)}/G$ of internal normal subgroupoids of $G$ and the category
$\NSCM/(A,B,\mu)$ of normal subcrossed modules of $(A,B,\mu)$ are equivalent.
\end{Cor}

We now define internal quotient groupoid as follows.

\begin{Def}\label{DefQuotientgroupgpd} {\em  Let $G$ be an internal groupoid in $\C$ and
$N$ an internal  normal subgroupoid of $G$. Let  $(A,B,\alpha)$ and
$(S,T,\sigma)$ be  respectively the crossed modules corresponding
to $G$ and $N$.  Then the internal groupoid corresponding to the
quotient crossed module    $(A/S,B/T,\rho)$ is called an {\em internal quotient
groupoid} and denoted by $G_N$.    }\end{Def}

So the set of objects and the set of morphisms of $G_N$ are respectively the
quotient groups with operations $G_0/N_0$ and $G_1/N_1$.

Here note that in the internal quotient groupoid $G_N$, the internal normal subgroupoid $N$ is
not wide and so is not a normal subgroupoid in the sense of Higgins \cite{Hi} and Brown
\cite{Br1}. If $N$ is wide then $G_N$ becomes a singular object, i.e. $G_N$ is an abelian group and
$[g]\star [g_1]=0$ for all $[g],[g_1]\in G_N$, in $\C$. Therefore the internal quotient
groupoid $G_N$ defined in Definition \ref{DefQuotientgroupgpd}  is not a quotient groupoid $G/N$  in the sense of Higgins \cite{Hi}.

In the following theorem we compare and relate quotient groupoid and  internal quotient groupoid.

\begin{Theo}
Let $G$ be an internal groupoid in $\C$ and $N$ an internal normal subgroupoid of $G$.
If $N$ is wide in $G$, then $N$ is a normal subgroupoid of $G$ and also
the quotient groupoid $G/N$ and the internal quotient groupoid $G_N$ are same if and only
if $N$ is transitive.
\end{Theo}

\begin{Prf}
Let $g\in G(x,y)$ and $n\in N(x)$. Then
\begin{align*}
 g\circ n &= n-1_x+g \\
  &= n-1_x+1_y-1_y+g\\
  &= n'-1_y+g \\
  &= n'\circ g \in N(y)\circ g
\end{align*}
where $n'=n-1_x+1_y\in N(y)$ since $N_1$ is an ideal of $G_1$. Hence $g\circ N(x)\subset N(y)\circ g$ and similarly
$N(y)\circ g\subset g\circ N(x)$. Therefore if $N$ is wide, then it  is a normal subgroupoid of $G$.

Let the quotient groupoid $G/N$ and the internal quotient groupoid $G_N$ be  same. Since $N$ is wide in $G$
then $G_N$, and hence $G/N$, has only one object. Hence $N$ is transitive.

On the other hand there are two equivalence relations on the set of morphisms  $G_1$.

$g \sim _1 g'$ if and only if there exists $n \in N_1$ such that
$g=g'+n$. Thus $G_1/\sim_1$ gives the quotient by the group with operations structure.

$g \sim_2 g'$ if and only if there are arrows  $m,n \in N_1$ such that
$g=n \circ g' \circ m$. Thus $G_1/ \sim_2$ gives the groupoid quotient.

If   $N$  is transitive, then these two groupoids  have only one object, i.e, they are singular objects in $\C$.
Suppose that $g \sim_2 g'$ in $G_1$ so that there  are $m,n \in N_1$ such that
$g=n\circ g' \circ m$. Then also
\begin{align*}
g & =  n\circ g' \circ m \\
  & =  n - 1_{d_0(n)} + g' - 1_{d_0(g')} + m \\
  & =  g' - g' + n - 1_{d_0(n)} + g' - 1_{d_0(g')} + m \\
  & =  g' + (- g' + n - 1_{d_0(n)} + g' - 1_{d_0(g')} + m).
\end{align*}
Since $N$ is an internal normal subgroupoid of $G$, \[(- g' + n - 1_{d_0(n)} + g' - 1_{d_0(g')} + m)\in N_1\]
so writing $n_1$ for $(- g' + n - 1_{d_0(n)} + g' - 1_{d_0(g')} + m)$ it follows that $g = g'+n_1$ and hence $g \sim_1 g'$.

On the other hand, if $g \sim_1 g'$, then there is an $n \in N_1$ such
that $g = g'+n$. Suppose  $g\in G(y,v)$ and $n\in N(x,u)$. Then $g\in G(y+x,v+u)$.
Since $N$ is wide in $G$, then $1_x \in N_1$ for all $x \in G_0$.
Since  $N$ is transitive then for each  $x\in N_0$ we can choose a
morphism $0_x\in N(0,x)$. Hence
\[[(1_v+n)\circ (1_v+0_x)]\circ  g' \circ (1_y + 0_x^{-1})= g\]

This implies that  $g \sim_2 g'$. Hence the quotient groupoid $G/N$ and the internal quotient groupoid $G_N$ are the same.
\end{Prf}

The following result characterize the normal objects in the category of internal groupoids in $\C$.
\begin{Theo}
Let $G$ and $N$ be two internal groupoids in $\C$. Then $N$ is an internal normal subgroupoid
of $G$ if and only if there exist an internal groupoid morphism $f\colon G\rightarrow H$
such that $\Ker f=N$.
\end{Theo}

\begin{Prf}
We know that the kernel of a  morphism $f\colon G\rightarrow H$  internal groupoids is an internal normal subgroupoid of $G$.

Conversely if $N$ is an internal normal subgroupoid of $G$, then  we can  obtain  the internal quotient groupoid
$G_N$  and quotient morphism   $p\colon G\rightarrow G_N$ of internal groupoids in $\C$ such that  $\Ker p=N$. \end{Prf}
\section{Normal and quotient cat$^1$-groups with operations}

We now generalize the notion of cat$^{1}$-group to the groups  with operations as follows.

\begin{Def} {\em
Let $A$ be an object in $\C$ and $s,t$  two endomorphisms of $A$ in $\C$. If the following are satisfied for all $\star\in\Omega_2'$, then $(A,s,t)$ is called a \emph{cat$^{1}$-group with operations}
or \emph{cat$^{1}$-group object in} $\C$.
\begin{enumerate} [label={{ (\roman{*})} }, leftmargin=1cm]
\item $st=t$ , $ts=s$;
\item $[\Ker s, \Ker t]=0$;
\item $\Ker s \star \Ker t=0$.
\end{enumerate}
}
\end{Def}

\begin{Exam} {\em Here are some examples of cat$^{1}$-groups with operations:

\begin{enumerate} [label={{ (\roman{*})} }, leftmargin=1cm]

\item [(i)] An  object $A$ in $\C$ can be regarded as a cat$^{1}$-group with operations
by the endomorphisms  $s=1_A=t$.

\item [(ii)]  A singular object $A$ in $\C$, i.e., an abelian group such that
$a\star a'=0$ for all $a,a'\in A$ and $\star\in\Omega_2'$, is a cat$^{1}$-group
with operations where $s$ and $t$ are zero morphisms.

\item [(iii)]  If  $B$ is  an object of  $\C$, then  $A=B\times B$ is a cat$^{1}$-group with operations
by the endomorphisms  $s$ and $t$ defined by  $s(b,b_1)=(b,b)$ and $t(b,b_1)=(b_1,b_1)$.
\end{enumerate}
}
\end{Exam}

Let $(A,s,t)$ and $(A',s',t')$ be two cat$^{1}$-groups with operations. A morphism $f\colon A\rightarrow A'$
in $\C$ is called a {\em morphism } of cat$^{1}$-groups with operations if it is compatible with the endomorphisms
$s$ and $t$, i.e. if $fs=s'f$ and $ft=t'f$.

Hence  we can construct a category denoted $\Cat^{1}(\C)$ of cat$^{1}$-groups with operations and their  morphisms. We now generalize a result proved in cat$^1$-group case by Loday in \cite{Loday82} to the cat$^1$-groups with operations in $\C$.

\begin{Theo}\label{Eqofcrossedmodandcatgroup}
The category $\CMod(\C)$ of crossed modules and the category $\Cat^1(\C)$  of cat$^{1}$-groups
with operations in $\C$ are equivalent.
\end{Theo}

\begin{Prf}
We give a sketch proof  based on  that of group case given  in \cite{Loday82}. First of all define a functor
\[\delta\colon \Cat^1(\C)\longrightarrow \CMod(\C)\] as follows: For an object  $(A,s,t)$ of $\Cat^1(\C)$,
$\delta(A,s,t)$ is $(\Ker s, \Im s, t_{|\Ker s})$ which is an object of $\CMod(\C)$ where the actions of $\Im s$ on $\Ker s$
are
\begin{eqnarray*}
  s(a)\cdot a_1 &=& s(a)+a_1-s(a) \\
  s(a)\star a_1 &=& s(a)\star a_1
\end{eqnarray*}
for $a,a_1\in A$ and  $\star\in\Omega_2'$. If $f\colon(A,s,t)\rightarrow(A',s',t')$ is a morphism in $\Cat^1(\C)$,
then \[\delta(f)=(f_{|\Ker s},f_{|\Im s})\colon (\Ker s, \Im s, t_{|\Ker s})\rightarrow(\Ker s', \Im s', t'_{|\Ker s'})\]
is a morphism in $\CMod(\C)$.

Conversely define a functor \[\theta\colon \CMod(\C) \longrightarrow \Cat^1(\C)\] in the following way: Let $(A,B,\alpha)$
be an object of $\CMod(\C)$. Then $\theta(A,B,\alpha)=(A\rtimes B,s,t)$ where $s(a,b)=(0,b)$ and $t(a,b)=(0,\alpha(a)+b)$.
If $(f,g)\colon(A,B,\alpha)\rightarrow (A',B',\alpha')$ is a morphism in $\CMod(\C)$, then
$f\times g\colon(A\rtimes B,s,t)\rightarrow (A'\rtimes B',s',t')$ is a morphism in $\Cat^1(\C)$.

Moreover we obtain a natural equivalence $\eta\colon 1_{\CMod(\C)}\rightarrow \delta\theta$ , where if $(A,B,\alpha)$
is a crossed module in $\C$, then $\eta_{(A,B,\alpha)}$ is given   by $a\mapsto (a,0)$ on $A$ and by $b\mapsto (0,b)$  on $B$.
Another natural equivalence $\mu\colon \theta\delta\rightarrow 1_{\Cat^1(\C)}$ is defined as follows: Let $(A,s,t)$ be
an object of $\Cat^1(\C)$. A morphism $\mu_{(A,s,t)}\colon \theta\delta((A,s,t))\rightarrow (A,s,t)$ is given by
$(a,s(b))\mapsto a+s(b)$. It is easy to verify that this morphism is an isomorphism in $\Cat^1(\C)$. The rest of
the proof is straightforward.
\end{Prf}

As a result of  Theorem \ref{Theocatequivalencess} and Theorem \ref{Eqofcrossedmodandcatgroup} we obtain the following corollary.

\begin{Cor} \label{Eqcatgroupandinternalcat}
The category  $\Cat^1(\C)$  of cat$^{1}$-groups in $\C$   and the category $\Cat(\C)$  of internal groupoids in $\C$ are equivalent.
\end{Cor}

By Theorem \ref{Eqofcrossedmodandcatgroup} we relate the subobjects in these categories as follows:

\begin{Theo}
Let $(S,T,\sigma)$ be a subcrossed module of a crossed module $(A,B,\alpha)$ in $\C$.  Suppose that  $(H,s_H,t_H)$ and $(G,s_G,t_G)$ are respectively
the cat$^1$-groups with operations corresponding to these crossed modules. Then $H$ is a subobject of $G$; and  $s_H$ and $t_H$ are respectively   the restrictions
of $s_G$ and $t_G$.
\end{Theo}

\begin{Prf} By the detailed  proof of Theorem \ref{Eqofcrossedmodandcatgroup}, we know that $H=S\rtimes T$ and $G=A\rtimes B$; and
 $s_H$ and $t_H$ are respectively the restrictions of $s_G$ and $t_G$.
Since  $(S,T,\sigma)$ is a subcrossed module of $(A,B,\alpha)$ in $\C$, it follows that $S\rtimes T$ is a subobject of $A\rtimes B$ in $\C$.
\end{Prf}

Hence we can state the notion of subobjects in $\Cat^{1}(\C)$ as follows:

\begin{Def} {\em
Let $(G,s_G,t_G)$ and $(H,s_H,t_H)$ be two objects of $\Cat^{1}(\C)$. If $H$ is a subobject of $G$ in $\C$;
and $s_H$ and $t_H$ are respectively the restrictions of $s_G$ and $t_G$ to $H$, then $(H,s_H,t_H)$ is called a
\emph{subcat$^1$-group with operations} or {\em subobject} of $(G,s_G,t_G)$ in $\Cat^{1}(\C)$.}\end{Def}

Similarly by Theorem \ref{Eqofcrossedmodandcatgroup} we relate the normal subobjects in these categories as follows:

\begin{Theo} \label{Thenormalsubcrosmod}
Let $(S,T,\sigma)$ be a normal subcrossed module of a crossed module $(A,B,\alpha)$ in $\C$.  Suppose that  $(N,s_N,t_N)$
and $(G,s_G,t_G)$ are respectively the cat$^1$-groups with operations corresponding to these crossed modules. Then $(N,s_N,t_N)$
is an  ideal of $G$.
\end{Theo}

\begin{Prf} By the proof of Theorem \ref{Eqofcrossedmodandcatgroup}, we know that $N=S\rtimes T$ and $G=A\rtimes B$ and
$N$ is an ideal of $G$ by the Proposition \ref{exactseqofcrosmod}.
\end{Prf}

Hence we can state the notion of normal subobject in $\Cat^{1}(\C)$ as follows:

\begin{Def}{\em
Let $(G,s_G,t_G)$ and $(N,s_N,t_N)$ be two objects in $\Cat^{1}(\C)$. If $N$ is an ideal of $G$ in $\C$;
and $s_N$ and $t_N$ are respectively the restrictions of $s_G$ and $t_G$ to $N$ , then $(N,s_N,t_N)$ is called a
\emph{normal subcat$^1$-group with operations} or {\em normal subobject} of $(G,s,t)$ in $\Cat^{1}(\C)$.}\end{Def}

For example if   $f\colon (G,s_G,t_G)\rightarrow(H,s_H,t_H)$ is  a morphism of cat$^1$-groups with operations in $\Cat^1(\C)$, then  the kernel \[\Ker f=\{g\in G ~|~ f(g)=0_H\in G\}\] of $f$ along with the endomorphisms $s_{\mid \Ker f}$ and $t_{\mid \Ker f}$ is a normal object of $(G,s_G,t_G)$.

\begin{Theo}
Let $(G,s_G,t_G)$ be an object in $\Cat^{1}(\C)$ and $(N,s_N,t_N)$ be a normal subcat$^1$-group
with operations of $(G,s_G,t_G)$. If $(A,B,\alpha)$ and $(S,T,\sigma)$ are  respectively  the  crossed
modules corresponding to $(G,s_G,t_G)$ and $(N,s_N,t_N)$. Then $(S,T,\sigma)$ is a normal subcrossed
module of $(A,B,\alpha)$ in $\C$.
\end{Theo}

\begin{Prf}
By the proof of Theorem \ref{Eqofcrossedmodandcatgroup} we have that $A=\Ker s_G$, $B=\Im s_G$, $S=\Ker s_N=A\cap N$
and $T=\Im s_N=s(N)$. We need to show that $(S,T,\sigma)$  satisfies the conditions of Definition
\ref{nrmlsbcrsmd}.
\begin{enumerate}[label={NCM\arabic{*}.}, leftmargin=2cm]
\item If  $g\in G$ and $n\in N$, then $s(g)+s(n)-s(g) = s(g+n-g)\in s(N)=T$ and $s(g)\star s(n)=s(g\star n)\in s(N)=T$
since $N$ is an ideal of $G$. Hence  $T$ is an ideal of $B$.
\item If $s(g)\in B$ and $n\in S$, then  by the proof of Theorem \ref{Eqofcrossedmodandcatgroup}
    $s(g) \cdot n  = s(g)+n-s(g)\in N$ since $N$ is an ideal of $G$.  Further
    \begin{align*}
        s(s(g)\cdot n) & = s(s(g)+n-s(g)) \\
                       & = s(s(g))+s(n)-s(s(g))\\
                       & = s(g)+0-s(g)\\
                       & = 0.
    \end{align*}
    and therefore  $s(g) \cdot n \in S$.
\item For  $s(n)\in T$ and $g\in A$ we have that
    \begin{align*}
        s((s(n)\cdot g)-g) & = s(s(n)+g-s(n)-g) \\
                           & = s(s(n))+s(g)-s(s(n))-s(g)\\
                           & = s(n)+0-s(n)-0 \\
                           & = 0
    \end{align*}
    and so $(s(n)\cdot g)-g\in A$. Since $g\in G$, $s(n)\in N$ and $N$ is an
    ideal of $G$ these imply  $g-s(n)-g\in N$ and since $s(n)\in N$ it implies that
    $s(n)+g-s(n)-g=(s(n)\cdot g)-g\in N$.  Hence $(s(n)\cdot g)-g\in S$.
\item Let $s(g)\in B$ and $n\in S$. Then by the proof of Theorem \ref{Eqofcrossedmodandcatgroup}
    $s(g) \star n  = s(g)\star n\in N$ since $N$ is an ideal of $G$. Further
    \begin{align*}
        s(s(g)\star n) & = s(s(g)\star n) \\
                       & = s(s(g))\star s(n)\\
                       & = s(g)\star 0\\
                       & = 0
    \end{align*}
    and thus $s(g) \star n \in S$.
\item For $s(n)\in T$ and $g\in A$ we have
    \begin{align*}
        s(s(n)\star g) & = s(s(n)\star g) \\
                           & = s(s(n))\star s(g)\\
                           & = s(n)\star 0 \\
                           & = 0
    \end{align*}
    and so $s(n)\star g\in A$. Since $g\in G$, $s(n)\in N$ and $N$ is an
    ideal of $G$ it follows that $s(n)\star g\in N$. Hence  $s(n)\star g\in S$.
\end{enumerate}
Therefore $(S,T,\sigma)$ is a normal subcrossed module of $(A,B,\alpha)$.
\end{Prf}

We now  can construct the quotient objects in the category $\Cat^{1}(\C)$ as follows:
Let $(G,s_G,t_G)$ be an object in $\Cat^{1}(\C)$ and $(N,s_N,t_N)$ a normal subobject of $(G,s_G,t_G)$. Then  $(G/N,s_{G/N},t_{G/N})$  becomes a
cat$^1$-group with operations by  the induced endomorphisms in $\C$
    \begin{align*}
        s_{G/N}([g]) & = [s_G(g)] \\
        t_{G/N}([g]) & = [t_G(g)].
    \end{align*}
This cat$^1$-group with operations is called the \emph{quotient cat$^1$-group with operations} of $(G,s_G,t_G)$
by $(N,s_N,t_N)$.

In the following theorem we prove that normal subcat$^1$-groups with operations are really the normal objects in the category  $\Cat^1(\C)$.
\begin{Theo}
Let $(G,s_G,t_G)$ and $(N,s_N,t_N)$ be two cat$^1$-groups with operations. Then $(N,s_N,t_N)$ is a normal
subcat$^1$-group with operations of $(G,s_G,t_G)$ if and only if it is a kernel of some morphism
 $f\colon (G,s_G,t_G)\rightarrow(H,s_H,t_H)$ of cat$^1$-groups with operations.
\end{Theo}

\begin{Prf}
We know the kernel of a morphism $f\colon (G,s_G,t_G)\rightarrow(H,s_H,t_H)$ in $\Cat^1(\C)$ is a normal subcat$^1$-group with
operations of $(G,s_G,t_G)$.

Conversely if  $(N,s_N,t_N)$ is a normal subcat$^1$-group with operations of $(G,s_G,t_G)$, then $p\colon (G,s_G,t_G)\rightarrow (G/N,s_{G/N},t_{G/N})$ is a morphism in $\Cat^1(\C)$ and $\Ker p=(N,s_N,t_N)$.
\end{Prf}

\section{Covering morphisms of internal groupoids and cat$^1$-groups with operations in $\C$}

Let $\mathbb{C}$ be an arbitrary category with pullbacks. The notion of {\em covering morphism} in $\mathbb{C}$  is defined in   \cite[pp.145]{Br-Da-Ha} as  a morphism  $p\colon \widetilde{G}\rightarrow G$ of internal
groupoids in $\mathbb{C}$ such that \[(p_1,d_0)\colon \widetilde{G}_1\longrightarrow G_1{_{d_0}\times_{p_0}}\widetilde{G}_0\]
is an isomorphism in $\mathbb{C}$. Equivalently it  is defined in  \cite{Ak-Al-Mu-Sa} as   an internal groupoid morphism $p\colon \widetilde{G}\rightarrow G$ in  $\mathbb{C}$ such that  $p$ is a covering morphism on the underlying groupoids.

Using the equivalence of the categories in Theorem \ref{Theocatequivalencess},  in \cite{Ak-Al-Mu-Sa}  a {\em covering morphism} of crossed modules is defined as a morphism  $(f,g)\colon (\widetilde{A},\widetilde{B},\tilde{\alpha})\rightarrow(A,B,\alpha)$ of crossed modules in $\C$ such that
if $f\colon \widetilde{A}\rightarrow A$ is an isomorphism in $\C$.

We now  give a result for coverings of  internal quotient groupoids.

\begin{Prop}
Let $p\colon \widetilde{G}\rightarrow G$ be a covering morphism in $\C$, $N$ an internal normal subgroupoid of $G$ and $\widetilde{N}=p^{-1}(N)$.
Then $\widetilde{N}$ is an internal normal subgroupoid of $\widetilde{G}$ and the induced morphism $p_{\ast}\colon \widetilde{G}_{\widetilde{N}}\rightarrow G_{N}$ is a covering morphism in $\C$.
\end{Prop}

\begin{Prf}
It is straightforward to prove that $\widetilde{N}=p^{-1}(N)$ is an internal normal subgroupoid of $\widetilde{G}$. Let
$(f,g)\colon (\widetilde{A},\widetilde{B},\tilde{\alpha})\rightarrow (A,B,\alpha)$ be the morphism of crossed modules corresponding
to $p$. Since $p$ is a covering morphism of internal groupoids in $\C$ then  $f\colon \widetilde{A}\rightarrow A$
is an isomorphism in $\C$.  Therefore $f_{\star}\colon \widetilde{A}/\widetilde{S}\rightarrow A/S$ is an isomorphism in $\C$
where $\widetilde{S}=\Ker d_{0_{|\widetilde{N}_1}}$ and $S=\Ker d_{0_{|N_1}}$. Hence $p_{\ast}\colon \widetilde{G}_{\widetilde{N}}\rightarrow G_N$
is a covering morphism in $\C$.
\end{Prf}

 Let  $G$  an internal groupoid in $\C$ and $(A,B,\alpha)$ be the crossed module corresponding to $G$ by Theorem \ref{Theocatequivalencess}. Then it is proved in \cite{Ak-Al-Mu-Sa} that the category  $\Cov_{\Cat(\C)}/G$ of covering morphisms based on $G$ and the category $\Cov_{\CMod(\C)}/(A,B,\alpha)$ of covering morphisms of crossed module based on $(A,B,\alpha)$ are equivalent.

 We now give a parallel result for cat$^1$-groups with operations.

\begin{Def}{\em A morphism $p\colon (\widetilde{G},s_{\widetilde{G}},t_{\widetilde{G}})\rightarrow (G,s_G,t_G)$ of cat$^1$-groups with operations  is called a \emph{covering morphism}  in $\Cat^1(\C)$ if the restriction $p_{|\Ker s_{\widetilde{G}}}\colon \Ker s_{\widetilde{G}}\rightarrow \Ker s_{G}$  of $p$ to  $\Ker s_{\widetilde{G}}$ is an isomorphism in $\C$.
}\end{Def}

Let   $\Cov_{\Cat^1(\C)}/(G,s_G,t_G)$ be the  category of covering morphisms   in $\Cat^1(\C)$ based on a cat$^1$-group with operations $(G,s_G,t_G)$.

We finally give the following Theorem.

\begin{Theo} Let $(A,B,\alpha)$ be a crossed module and let  $(G,s_G,t_G)$  be the  corresponding cat$^1$-group with operations. Then the category
$\Cov_{\CMod(\C)}/(A,B,\alpha)$ of covers of $(A,B,\alpha)$ and the category $\Cov_{\Cat^1(\C)}/(G,s_G,t_G)$
are equivalent.
\end{Theo}

\begin{Prf} A  functor $\eta\colon \Cov_{\CMod(\C)}/(A,B,\alpha)\rightarrow \Cov_{\Cat^1(\C)}/(G,s_G,t_G)$ is defined as follows:
Let $(f,g)\colon (\widetilde{A},\widetilde{B},\tilde{\alpha})\rightarrow(A,B,\alpha)$ be an object in $\Cov_{\CMod(\C)}/(A,B,\alpha)$. Then
$\eta(f,g)=f\times g\colon \widetilde{A}\rtimes \widetilde{B}\rightarrow A\rtimes B$,  $s_{\widetilde{A}\rtimes \widetilde{B}}(\tilde{a},\tilde{b})=(0,\tilde{b})$ and
$s_{A\rtimes B}(a,b)=(0,b)$. Here $\Ker s_{\widetilde{A}\rtimes \widetilde{B}}=\widetilde{A}\times\{0\}$, $\Ker s_{A\rtimes B}=A\times\{0\}$ and
the restriction of $f\times g$ to $\Ker s_{\widetilde{A}\rtimes \widetilde{B}}$ is $f\times 0$. Since $(f,g)$ is a covering morphism then
$f$ and hence $f\times 0$ is an isomorphism in $\C$. Thus $\eta(f,g)=f\times g$ becomes an object in $\Cov_{\Cat^1(\C)}/(G,s_G,t_G)$.

Conversely a  functor $\theta\colon \Cov_{\Cat^1(\C)}/(G,s_G,t_G) \rightarrow \Cov_{\CMod(\C)}/(A,B,\alpha)$ is
defined in the following way: Let $p\colon (\widetilde{G},s_{\widetilde{G}},t_{\widetilde{G}})\rightarrow (G,s_G,t_G)$ be an object in
$\Cov_{\Cat^1(\C)}/(G,s_G,t_G)$. Then $\theta(p)=(p_{|\Ker s_{\widetilde{G}}},p_{|\Im s_{\widetilde{G}}})$ and clearly $\theta(p)$
is an isomorphism and hence an object in $\Cov_{\CMod(\C)}/(A,B,\alpha)$.

The other details of the equivalence of the categories can be checked.
\end{Prf}

\end{document}